# ASYMPTOTICS IN RANDOMIZED URN MODELS

By Zhi-Dong Bai[1] and Feifang Hu[2]

*Northeast Normal University and National University of Singapore, and University of Virginia*

This paper studies a very general urn model stimulated by designs in clinical trials, where the number of balls of different types added to the urn at trial $n$ depends on a random outcome directed by the composition at trials $1, 2, \ldots, n-1$. Patient treatments are allocated according to types of balls. We establish the strong consistency and asymptotic normality for both the urn composition and the patient allocation under general assumptions on random generating matrices which determine how balls are added to the urn. Also we obtain explicit forms of the asymptotic variance–covariance matrices of both the urn composition and the patient allocation. The conditions on the nonhomogeneity of generating matrices are mild and widely satisfied in applications. Several applications are also discussed.

**1. Introduction.** In designing a clinical trial, the limiting behavior of the patient allocation to several treatments during the process is of primary consideration. Suppose patients arrive sequentially from a population. Adaptive designs in clinical trials are inclining to assign more patients to better treatments, while seeking to maintain randomness as a basis for statistical inference. Thus the cumulative information of the responses of treatments on previous patients will be used to adjust treatment assignment to coming patients. For this purpose, various urn models [Johnson and Kotz (1977)] have been proposed and used extensively in adaptive designs [for more references, see Zelen (1969), Wei (1979), Flournoy and Rosenberger (1995) and Rosenberger (1996)].

One large family of randomized adaptive designs is based on the generalized Friedman's urn (GFU) model [Athreya and Karlin (1967, 1968),

Received June 2003; revised March 2004.

[1] Supported by NSFC Grant 201471000 and NUS Grant R-155-000-030-112.

[2] Supported by NSF Grant DMS-0204232 and NUS Grant R-155-000-030-112.

also called the generalized Pólya urn (GPU) in the literature]. The model can be described as follows. Consider an urn containing balls of $K$ types, respectively, representing $K$ "treatments" in a clinical trial. These treatments are to be assigned sequentially in $n$ stages. At the beginning, the urn contains $\mathbf{Y}_0 = (Y_{01}, \ldots, Y_{0K})$ balls, where $Y_{0k}$ denotes the number of balls of type $k$, $k = 1, \ldots, K$. At stage $i$, $i = 1, \ldots, n$, a ball is randomly drawn from the urn and then replaced. If the ball is of type $q$, then the treatment $q$ is assigned to the $i$th patient, $q = 1, \ldots, K$, $i = 1, \ldots, n$. We then wait until we observe a random variable $\xi(i)$, which may include the response and/or other covariates of patient $i$. After that, an additional $D_{qk}(i)$ balls of type $k$, $k = 1, \ldots, K$, are added to the urn, where $D_{qk}(i)$ is some function of $\xi(i)$. This procedure is repeated throughout the $n$ stages. After $n$ splits and generations, the urn composition is denoted by the row vector $\mathbf{Y}_n = (Y_{n1}, \ldots, Y_{nK})$, where $Y_{nk}$ represents the number of balls of type $k$ in the urn after the $n$th split. This relation can be written as the following recursive formula:

$$\mathbf{Y}_n = \mathbf{Y}_{n-1} + \mathbf{X}_n \mathbf{D}_n,$$

where $\mathbf{X}_n$ is the result of the $n$th draw, distributed according to the urn composition at the previous stage; that is, if the $n$th draw is a type-$k$ ball, then the $k$th component of $\mathbf{X}_n$ is 1 and other components are 0. Furthermore, write $\mathbf{N}_n = (N_{n1}, \ldots, N_{nK})$, where $N_{nk}$ is the number of times a type-$k$ ball was drawn in the first $n$ stages, or equivalently, the number of patients who receive the treatment $k$ in the first $n$ patients.

For notation, let $\mathbf{D}_i = \langle\langle D_{qk}(i), q, k = 1, \ldots, K \rangle\rangle$ and let $\mathcal{F}_i$ be the sequence of increasing $\sigma$-fields generated by $\{\mathbf{Y}_j\}_{j=0}^i$, $\{\mathbf{X}_j\}_{j=1}^i$ and $\{\mathbf{D}_j\}_{j=1}^i$. Define $\mathbf{H}_i = \langle\langle E(D_{qk}(i)|\mathcal{F}_{i-1}), q, k = 1, \ldots, K \rangle\rangle$, $i = 1, \ldots, n$. The matrices $\mathbf{D}_i$ are called *addition rules* and $\mathbf{H}_i$ *generating matrices*. In practice, the addition rule $\mathbf{D}_i$ often depends only on the treatment on the $i$th patient and its outcome. In these cases, the addition rules $\mathbf{D}_i$ are i.i.d. (independent and identically distributed) and the generating matrices $\mathbf{H}_i = \mathbf{H} = E\mathbf{D}_i$ are identical and nonrandom. But in some applications, the addition rule $\mathbf{D}_i$ depends on the total history of previous trials [see Andersen, Faries and Tamura (1994) and Bai, Hu and Shen (2002)]; then the general generating matrix $\mathbf{H}_i$ is the conditional expectation of $\mathbf{D}_i$ given $\mathcal{F}_{i-1}$. Therefore, the general generating matrices $\{\mathbf{H}_i\}$ are usually random. In this paper, we consider this general case. Examples are considered in Section 5.

A GFU model is said to be *homogeneous* if $\mathbf{H}_i = \mathbf{H}$ for all $i = 1, 2, 3, \ldots$. In the literature, research is focused on asymptotic properties of $\mathbf{Y}_n$ for homogeneous GFU. First-order asymptotics for homogeneous GFU models are determined by the generating matrices $\mathbf{H}$. In most cases, $\mathbf{H}$ is an irreducible nonnegative matrix, for which the maximum eigenvalue is unique and positive (called the maximal eigenvalue in the literature) and its corresponding left eigenvector has positive components. In some cases, the entries of $\mathbf{H}$ may not be all nonnegative (e.g., when there is no replacement after the draw), and we may assume that the matrix $\mathbf{H}$ has a unique maximal eigenvalue $\lambda$ with associated left eigenvector $\mathbf{v} = (v_1, \ldots, v_K)$ with $\sum v_i = 1$. Under the following assumptions:

(i) $\Pr\{D_{qk} = 0, k = 1, \ldots, K\} = 0$ for every $q = 1, \ldots, K$,
(ii) $D_{qk} \geq 0$ for all $q, k = 1, \ldots, K$,
(iii) $\mathbf{H}$ is irreducible,

Athreya and Karlin (1967, 1968) prove that

$$(1.1) \qquad \frac{N_{nk}}{n} \to v_k \quad \text{and} \quad \frac{Y_{nk}}{\sum_{q=1}^{K} Y_{nq}} \to v_k$$

almost surely as $n \to \infty$.

Let $\lambda_1$ be the eigenvalue with a second largest real part, associated with a right eigenvector $\xi$. If $\lambda > 2\operatorname{Re}(\lambda_1)$, Athreya and Karlin (1968) show that

$$(1.2) \qquad n^{-1/2}\mathbf{Y}_n\xi' \to N(0, c)$$

in distribution, where $c$ is a constant. When $\lambda = 2\operatorname{Re}(\lambda_1)$ and $\lambda_1$ is simple, then (1.2) holds when $n^{-1/2}$ is replaced by $1/\sqrt{n\ln(n)}$. Asymptotic results under various addition schemes are considered in Freedman (1965), Mahmoud and Smythe (1991), Holst (1979) and Gouet (1993).

Homogeneity of the generating matrix is often not the case in clinical trials, where patients may exhibit a drift in characteristics over time. Examples are given in Altman and Royston (1988), Coad (1991) and Hu and Rosenberger (2000). Bai and Hu (1999) establish the weak consistency and the asymptotic normality of $\mathbf{Y}_n$ under GFU models with nonhomogeneous generating matrices $\mathbf{H}_i$. [In that paper, it is assumed that $\mathbf{H}_i = E\mathbf{D}_i$, so $\mathbf{H}_i$ are fixed (not random) matrices.] They consider the following GFU model (GFU1): $\sum_{k=1}^{K} D_{qk}(i) = c_1 > 0$, for all $q = 1, \ldots, K$ and $i = 1, \ldots, n$, the total number of balls added at each stage is a positive constant. They assume there is a nonnegative matrix $\mathbf{H}$ such that

$$(1.3) \qquad \sum_{i=1}^{\infty} \frac{\alpha_i}{i} < \infty,$$

where $\alpha_i = \|\mathbf{H}_i - \mathbf{H}\|_\infty$.

In clinical trials, $N_{nk}$ represents the number of patients assigned to the treatment $k$ in the first $n$ trials. Doubtless, the asymptotic distribution and asymptotic variance of $\mathbf{N}_n = (N_{n1}, \ldots, N_{nK})$ is of more practical interest than the urn compositions to sequential design researchers. As Athreya and Karlin [(1967), page 275] said, "It is suggestive to conjecture that $(N_{n1}, \ldots, N_{nK})$ properly normalized is asymptotically normal. This problem is open." The problem has stayed open for decades due to mathematical complexity. One of our main goals of this paper is to present a solution to this problem.

Smythe (1996) defined the extended Pólya urn (EPU) (homogeneous) models, satisfying $\sum_{k=1}^{K} E(D_{qk}) = c_1 > 0, q = 1, \ldots, K$; that is, the expected total number of balls added to the urn at each stage is a positive constant. For EPU models, Smythe (1996) established the weak consistency and the asymptotic normality of $\mathbf{Y}_n$ and $\mathbf{N}_n$ under the assumptions that the eigenvalues of the generating matrix $\mathbf{H}$ are simple. The asymptotic variance of $\mathbf{N}_n$ is a more important and difficult proposition [Rosenberger (2002)]. Recently, Hu and Rosenberger (2003) obtained an explicit relationship between the power and the variance of $\mathbf{N}_n$ in an adaptive design. To compare the randomized urn models with other adaptive designs, one just has to calculate and compare their variances. Matthews and Rosenberger (1997) obtained the formula for asymptotic variance for the randomized play-the-winner rule ($K = 2$) which was initially proposed by Wei and Durham (1978). A general formula for asymptotic variance of $\mathbf{N}_n$ was still an open problem [Rosenberger (2002)].

In this paper, we

(i) show the asymptotic normality of $\mathbf{N}_n$ for general $\mathbf{H}$;

(ii) obtain a general and explicit formula for the asymptotic variance of $\mathbf{N}_n$;

(iii) show the strong consistency of both $\mathbf{Y}_n$ and $\mathbf{N}_n$; and

(iv) extend these results to nonhomogeneous urn model with random generating matrices $\mathbf{H}_i$.

The paper is organized as follows. The strong consistency of $\mathbf{Y}_n$ and $\mathbf{N}_n$ is proved in Section 2 for both homogeneous and nonhomogeneous EPU models. Note that the GFU1 is a special case of EPU. The asymptotic normality of $\mathbf{Y}_n$ for homogeneous and nonhomogeneous EPU models is shown in Section 3 under the assumption (1.3). We consider cases where the generating matrix $\mathbf{H}$ has a general Jordan form. In Section 4, we consider the asymptotic normality of $\mathbf{N}_n = (N_{n1}, \ldots, N_{nK})$ for both homogeneous and nonhomogeneous EPU models. Further, we obtain a general and explicit formula for the asymptotic variance of $\mathbf{N}_n$.

The condition (1.3) in a nonhomogeneous urn model is widely satisfied in applications. In some applications [e.g., Bai, Hu and Shen (2002)], the generating matrix $\mathbf{H}_i$ may be estimates of some unknown parameters updated at each stage, for example, $\hat{\mathbf{H}}_i$ at $i$th stage. In these cases, we usually have $\alpha_i = O(i^{-1/2})$ in probability or $O(i^{-1/4})$ almost surely, so the condition (1.3) is satisfied. Also (1.3) is satisfied for the case of Hu and Rosenberger (2000). Some other applications are considered in Section 5.

**2. Strong consistency of $\mathbf{Y}_n$ and $\mathbf{N}_n$.** Using the notation defined in the Introduction, $\mathbf{Y}_n$ is a sequence of random $K$-vectors of nonnegative elements which are adaptive with respect to $\{\mathcal{F}_n\}$, satisfying

$$(2.1) \qquad E(\mathbf{Y}_i|\mathcal{F}_{i-1}) = \mathbf{Y}_{i-1}\mathbf{M}_i,$$

where $\mathbf{M}_i = \mathbf{I} + a_{i-1}^{-1}\mathbf{H}_i$, $\mathbf{H}_i = E(\mathbf{D}_i|\mathcal{F}_{i-1})$ and $a_i = \sum_{j=1}^{K} Y_{ij}$. Without loss of generality, we assume $a_0 = 1$ in the following study.

In the sequel, we need the following assumptions.

ASSUMPTION 2.1. The generating matrix $\mathbf{H}_i$ satisfies

$$(2.2) \qquad \begin{aligned} H_{qk}(i) &\geq 0 \qquad &&\text{for all } k, q \quad \text{and} \\ \sum_{k=1}^{K} H_{qk}(i) &= c_1 \qquad &&\text{for all } q = 1, \ldots, K, \end{aligned}$$

almost surely, where $H_{qk}(i)$ is the $(q, k)$-entry of the matrix $\mathbf{H}_i$ and $c_1$ is a positive constant. Without loss of generality, we assume $c_1 = 1$ throughout this work.

ASSUMPTION 2.2. The addition rule $\mathbf{D}_i$ is conditionally independent of the drawing procedure $\mathbf{X}_i$ given $\mathcal{F}_{i-1}$ and satisfies

$$(2.3) \quad E(D_{qk}^{2+\delta}(i)|\mathcal{F}_{i-1}) \leq C < \infty \qquad \text{for all } q, k = 1, \ldots, K \text{ and some } \delta > 0.$$

Also we assume that

$$(2.4) \quad \text{cov}[(D_{qk}(i), D_{ql}(i))|\mathcal{F}_{i-1}] \to d_{qkl} \qquad \text{for all } q, k, l = 1, \ldots, K,$$

where $\mathbf{d}_q = (d_{qkl})_{k,l=1}^{K}$, $q = 1, \ldots, K$, are some $K \times K$ positive definite matrices.

REMARK 2.1. Assumption 2.1 defines the EPU model [Smythe (1996)]; it ensures that the number of expected balls added at each stage is a positive constant. So after $n$ stages, the total number of balls, $a_n$, in the urn should be very close to $n$ ($a_n/n$ converges to 1).

The elements of the addition rule are allowed to take negative values in the literature, which corresponds to the situation of withdrawing balls from the urn. But, to avoid the dilemma that there are no balls to withdraw, only diagonal elements of $\mathbf{D}_i$ are allowed to take negative values, which corresponds to the case of drawing without replacement.

To investigate the limiting properties of $\mathbf{Y}_n$, we first derive a decomposition. From (2.1), it is easy to see that

$$(2.5) \qquad \begin{aligned} \mathbf{Y}_n &= (\mathbf{Y}_n - E(\mathbf{Y}_n|\mathcal{F}_{n-1})) + \mathbf{Y}_{n-1}\mathbf{M}_n \\ &= \mathbf{Q}_n + \mathbf{Y}_{n-1}\mathbf{G}_n + \mathbf{Y}_{n-1}(\mathbf{M}_n - \mathbf{G}_n) \\ &= \mathbf{Y}_0\mathbf{G}_1\mathbf{G}_2\cdots\mathbf{G}_n + \sum_{i=1}^{n} \mathbf{Q}_i\mathbf{B}_{n,i} + \sum_{i=1}^{n} \mathbf{Y}_{i-1}(\mathbf{M}_i - \mathbf{G}_i)\mathbf{B}_{n,i} \\ &= \mathbf{S}_1 + \mathbf{S}_2 + \mathbf{S}_3, \end{aligned}$$

where $\mathbf{Q}_i = \mathbf{Y}_i - E(\mathbf{Y}_i|\mathcal{F}_{i-1})$, $\mathbf{G}_i = \mathbf{I} + i^{-1}\mathbf{H}$ and $\mathbf{B}_{n,i} = \mathbf{G}_{i+1}\cdots\mathbf{G}_n$ with the convention that $\mathbf{B}_{n,n} = \mathbf{I}$ and $\mathcal{F}_0$ denotes the trivial $\sigma$-field.

We further decompose $\mathbf{S}_3$ as follows:

$$
\begin{aligned}
\mathbf{S}_3 &= \sum_{i=1}^{n} \mathbf{Y}_{i-1}(a_{i-1}^{-1}\mathbf{H}_i - i^{-1}\mathbf{H})\mathbf{B}_{n,i} \\
&= \sum_{i=1}^{n} a_{i-1}^{-1}\mathbf{Y}_{i-1}(\mathbf{H}_i - \mathbf{H})\mathbf{B}_{n,i} + \sum_{i=1}^{n} \frac{\mathbf{Y}_{i-1}}{a_{i-1}}\frac{i - a_{i-1}}{i}\mathbf{H}\mathbf{B}_{n,i} \\
&= \mathbf{S}_{31} + \mathbf{S}_{32}.
\end{aligned}
\tag{2.6}
$$

To estimate the above terms in the expansion, we need some preliminary results. First, we evaluate the convergence rate of $a_n$. To this end, we have the following theorem.

THEOREM 2.1.  *Under Assumptions* 2.1 *and* 2.2, (a) $a_n/n \to 1$ *a.s. as* $n \to \infty$, *and* (b) $n^{-\kappa}(a_n - n) \to 0$ *a.s. for any* $\kappa > 1/2$.

PROOF.   Let $e_i = a_i - a_{i-1}$ for $i \geq 1$. By definition, we have $e_i = \mathbf{X}_i\mathbf{D}_i\mathbf{1}$, where $\mathbf{X}_i$ is the result of the $i$th draw, multinomially distributed according to the urn composition at the previous stages; that is, the conditional probability that the $i$th draw is a ball of type $k$ (the $k$th component of $\mathbf{X}_i$ is 1 and other components are 0) given previous status is $Y_{i-1,k}/a_{i-1}$.

From Assumptions 2.1 and 2.2, we have

$$
E(e_i|\mathcal{F}_{i-1}) = 1
\tag{2.7}
$$

and

$$
\begin{aligned}
E(e_i^2) &= E[E(e_i^2|\mathcal{F}_{i-1})] = E[E(\mathbf{1}'\mathbf{D}_i'\mathbf{X}_i'\mathbf{X}_i\mathbf{D}_i\mathbf{1}|\mathcal{F}_{i-1})] \\
&= \mathbf{1}'E[E(\mathbf{D}_i'\mathbf{X}_i'\mathbf{X}_i\mathbf{D}_i|\mathcal{F}_{i-1})]\mathbf{1} \\
&= \mathbf{1}'E[E(\mathbf{D}_i'\operatorname{diag}(a_{i-1}^{-1}\mathbf{Y}_{i-1})\mathbf{D}_i|\mathcal{F}_{i-1})]\mathbf{1} \\
&= \sum_{q=1}^{K}\sum_{k=1}^{K}\sum_{l=1}^{K} E[(a_{i-1}^{-1}Y_{i-1,q})E(D_{qk}(i)D_{ql}(i)|\mathcal{F}_{i-1})] \\
&\leq CK^2,
\end{aligned}
\tag{2.8}
$$

so that

$$
a_n - n = a_0 + \sum_{i=1}^{n}(e_i - 1) = 1 + \sum_{i=1}^{n}(e_i - E(e_i|\mathcal{F}_{i-1}))
\tag{2.9}
$$

forms a martingale sequence.

From Assumption 2.2 and $\kappa > 1/2$, we have

$$
\sum_{i=1}^{\infty} E\left(\left(\frac{e_i - 1}{i^\kappa}\right)^2\Big|\mathcal{F}_{i-1}\right) < \infty.
$$

By three series theorem for martingales, this implies that the series

$$\sum_{i=1}^{\infty} \frac{e_i - 1}{i^\kappa}$$

converges almost surely. Then, by Kronecker's lemma,

$$\frac{1}{n^\kappa} \sum_{i=1}^{n} (e_i - 1) \to 0$$

almost surely. This completes the proof for conclusion (b) of the theorem.

The conclusion (a) is a consequence of conclusion (b). The proof of Theorem 2.1 is then complete. $\square$

ASSUMPTION 2.3. Assume that (1.3) holds almost surely. Suppose that the limit generating matrix $\mathbf{H}$, $K \times K$, is irreducible.

This assumption guarantees that $\mathbf{H}$ has the Jordan form decomposition

$$\mathbf{T}^{-1}\mathbf{H}\mathbf{T} = \mathbf{J} = \begin{pmatrix} 1 & 0 & \cdots & 0 \\ 0 & \mathbf{J}_1 & \cdots & 0 \\ \cdots & \cdots & \cdots & \cdots \\ 0 & 0 & \cdots & \mathbf{J}_s \end{pmatrix} \qquad \text{with } \mathbf{J}_t = \begin{pmatrix} \lambda_t & 1 & 0 & \cdots & 0 \\ 0 & \lambda_t & 1 & \cdots & 0 \\ \vdots & \cdots & \ddots & \ddots & \vdots \\ 0 & 0 & \cdots & \lambda_t & 1 \\ 0 & 0 & 0 & \cdots & \lambda_t \end{pmatrix},$$

where 1 is the unique maximal eigenvalue of the matrix $\mathbf{H}$. Denote the order of $\mathbf{J}_t$ by $\nu_t$ and $\tau = \max\{\mathrm{Re}(\lambda_1), \ldots, \mathrm{Re}(\lambda_s)\}$. We define $\nu = \max\{\nu_t : \mathrm{Re}(\lambda_t) = \tau\}$.

Moreover, the irreducibility of $\mathbf{H}$ also guarantees that the elements of the left eigenvector $\mathbf{v} = (v_1, \ldots, v_p)$ associated with the positive maximal eigenvalue 1 are positive. Thus, we may normalize this vector to satisfy $\sum_{i=1}^{p} v_i = 1$.

REMARK 2.2. Condition (1.3) in Assumption 2.3 is very mild, just slightly stronger than $\alpha_i \to 0$, for example, if the nonhomogeneous generating matrix $\mathbf{H}_i$ converges to a generating matrix $\mathbf{H}$ with a rate of $\log^{-1-c} i$ for some $c > 0$.

What we consider here is the general case where the Jordan form of the generating matrix $\mathbf{H}$ is arbitrary, relaxing the constraint of a diagonal Jordan form as usually assumed in the literature [see Smythe (1996)].

In some conclusions, we need the convergence rate of $\mathbf{H}_i$ as described in the following assumption.

ASSUMPTION 2.4.

$$(2.10) \qquad \|\mathbf{H}_i - E\mathbf{H}_i\| = \begin{cases} O(i^{-1/2}), & \text{if } \tau \neq \frac{1}{2}, \\ O(i^{-1/2}\log^{-1/2}(i+1)), & \text{if } \tau = \frac{1}{2}, \end{cases}$$

where $\|(a_{ij})\| = \sqrt{\sum_{ij} Ea_{ij}^2}$, for any random matrix $(a_{ij})$.

A slightly stronger condition is

$$(2.11) \qquad \|\mathbf{H}_i - E\mathbf{H}_i\| = o(i^{-1/2}).$$

REMARK 2.3. This assumption is trivially true if $\mathbf{H}_i$ is nonrandom. It is also true when $\mathbf{H}_i$ is a continuously differentiable matrix function of status at stage $i$, such as $\mathbf{Y}_i$, $\mathbf{N}_i$ or the relative frequencies of the success, and so on. These are true in almost all practical situations.

For further studies, we define

$$V_n = \begin{cases} \sqrt{n}, & \text{if } \tau < 1/2, \\ \sqrt{n}\log^{\nu-1/2} n, & \text{if } \tau = 1/2, \\ n^{\tau}\log^{\nu-1} n, & \text{if } \tau > 1/2. \end{cases}$$

THEOREM 2.2. Under Assumptions 2.1–2.3, for some constant $M$,

$$(2.12) \qquad E\|\mathbf{Y}_n - E\mathbf{Y}_n\|^2 \leq MV_n^2.$$

From this, for any $\kappa > \tau \vee \frac{1}{2}$, we immediately obtain $n^{-\kappa}(\mathbf{Y}_n - E\mathbf{Y}_n) \to 0$, a.s., where $a \vee b = \max(a,b)$. Also, if $\kappa = 1$ or the condition (1.3) is strengthened to

$$(2.13) \qquad \sum_{i=1}^{\infty} \frac{\alpha_i}{\sqrt{i}} < \infty,$$

then $E\mathbf{Y}_n$ in the above conclusions can be replaced by $n\mathbf{v}$. This implies that $n^{-1}\mathbf{Y}_n$ almost surely converges to $\mathbf{v}$, the same limit of $n^{-1}E\mathbf{Y}_n$, as $n \to \infty$.

PROOF. Without loss of generality, we assume $a_0 = 1$ in the following study. For any random vector, write $\|\mathbf{Y}\| := \sqrt{E\mathbf{Y}\mathbf{Y}'}$. Define $\mathbf{y}_n = (y_{n,1}, \ldots, y_{n,K}) = \mathbf{Y}_n\mathbf{T}$. Then, (2.12) reduces to

$$(2.14) \qquad \|\mathbf{y}_n - E\mathbf{y}_n\| \leq MV_n.$$

In Theorem 2.1, we have proved that $\|a_n - n\|^2 \leq CK^2n$ [see (2.9) and (2.8)]. Noticing that $Ea_n = n+1$, the proof of (2.12) further reduces to showing that, for any $j = 2, \ldots, K$,

$$(2.15) \qquad \|y_{n,j} - Ey_{n,j}\| \leq MV_n.$$

We shall prove (2.15) by induction.

Suppose $n_0$ is an integer and $M$ a constant such that

$$(2.16) \qquad \sum_{i=n_0}^{\infty} \frac{\alpha_i}{i} < \varepsilon, \qquad \sum_{i=n_0}^{\infty} \frac{\log^\nu i}{i^{5/4}} < \varepsilon,$$

$$M = \frac{C_1 + C_2 + C_3 + C_4 + C_5 + (C_3 + 2C_5)M_0}{1 - 3\varepsilon},$$

where $\varepsilon < 1/4$ is a prechosen small positive number, $M_0 = \max_{n \le n_0}\{\|y_{n,j} - Ey_{n,j}\|/V_n\}$ and the constants $C$'s are absolute constants specified later.

We shall complete the proof by induction. Consider $m > n_0$ and assume that $\|\tilde{\mathbf{y}} - E\tilde{\mathbf{y}}_n\| \le MV_n$ for all $n_0 \le n < m$.

By (2.5) and (2.6), we have

$$(2.17) \begin{aligned} y_{m,j} &= \mathbf{y}_0 \widetilde{\mathbf{B}}_{m,0,j} + \sum_{i=1}^{m} \widetilde{\mathbf{Q}}_i \widetilde{\mathbf{B}}_{m,i,j} + \sum_{i=1}^{m} \frac{\mathbf{y}_{i-1}}{i} \mathbf{W}_i \widetilde{\mathbf{B}}_{m,i,j} \\ &\quad + \sum_{i=1}^{m} \left( \frac{i - a_{i-1}}{i} \right) \frac{\mathbf{y}_{i-1}}{a_{i-1}} \mathbf{J} \widetilde{\mathbf{B}}_{m,i,j} + \sum_{i=1}^{m} \left( \frac{i - a_{i-1}}{i} \right) \frac{\mathbf{y}_{i-1}}{a_{i-1}} \mathbf{W}_i \widetilde{\mathbf{B}}_{m,i,j}, \end{aligned}$$

where $\widetilde{\mathbf{Q}}_i = \mathbf{Q}_i \mathbf{T}$, $\mathbf{W}_i = \mathbf{T}^{-1}(\mathbf{H}_i - \mathbf{H})\mathbf{T}$ and

$$(2.18) \begin{aligned} \widetilde{\mathbf{B}}_{n,i} &= \mathbf{T}^{-1}\mathbf{B}_{n,i}\mathbf{T} = (\mathbf{I} + (i+1)^{-1}\mathbf{J})\cdots(\mathbf{I} + n^{-1}\mathbf{J}) \\ &= \begin{pmatrix} \prod_{j=i+1}^{n}\left(1 + \frac{1}{j}\right) & 0 & \cdots & 0 \\ 0 & \prod_{j=i+1}^{n}(\mathbf{I} + j^{-1}\mathbf{J}_1) & \cdots & 0 \\ \cdots & \cdots & \cdots & \cdots \\ 0 & 0 & \cdots & \prod_{j=i+1}^{n}(\mathbf{I} + j^{-1}\mathbf{J}_s) \end{pmatrix}, \end{aligned}$$

and $\widetilde{\mathbf{B}}_{m,i,j}$ is the $j$th column of the matrix $\widetilde{\mathbf{B}}_{m,i}$.

In the remainder of the proof of the theorem, we shall frequently use the elementary fact that

$$(2.19) \qquad \prod_{j=i+1}^{n}\left(1 + \frac{\lambda}{j}\right) = \left(\frac{n}{i}\right)^{\lambda}\psi(n,i,\lambda),$$

where $\psi(n,i,\lambda)$ is uniformly bounded (say $\le \psi$) and tends to 1 as $i \to \infty$. In the sequel, we use $\psi(n,i,\lambda)$ as a generic symbol, that is, it may take different values at different appearances and is uniformly bounded (by $\psi$, say) and tends to 1 as $i \to \infty$. Based on this estimation, one finds that the $(h, h + \ell)$-element of the block matrix $\prod_{i=j+2}^{n}(\mathbf{I} + i^{-1}\mathbf{J}_t)$ is asymptotically equivalent to

$$(2.20) \qquad \frac{1}{\ell!}\left(\frac{j}{n}\right)^{-\lambda_t}\log^\ell\left(\frac{n}{j}\right)\psi(n,j,\lambda_t),$$

where $\lambda_t$ is the eigenvalue of $\mathbf{J}_t$.

By (2.17) and triangular inequality, we have

$$
\begin{aligned}
\|y_{m,j} &- Ey_{m,j}\| \\
&\leq \|\mathbf{y}_0 \widetilde{\mathbf{B}}_{m,0,j}\| + \left\|\sum_{i=1}^{m} \widetilde{\mathbf{Q}}_i \widetilde{\mathbf{B}}_{m,i,j}\right\| \\
&\quad + \left\|\sum_{i=1}^{m}\left(\frac{i-a_{i-1}}{i}\right)\frac{\mathbf{y}_{i-1}}{a_{i-1}}\mathbf{J}\widetilde{\mathbf{B}}_{m,i,j}\right\| + \sum_{i=1}^{m}\left\|\left(\frac{i-a_{i-1}}{i}\right)\frac{\mathbf{y}_{i-1}}{a_{i-1}}\mathbf{W}_i\widetilde{\mathbf{B}}_{m,i,j}\right\| \\
&\quad + \sum_{i=1}^{m}\left\|\frac{\mathbf{y}_{i-1}\mathbf{W}_i - E\mathbf{y}_{i-1}\mathbf{W}_i}{i}\widetilde{\mathbf{B}}_{m,i,j}\right\|.
\end{aligned} \tag{2.21}
$$

Consider the case where $1 + \nu_1 + \cdots + \nu_{t-1} < j \leq 1 + \nu_1 + \cdots + \nu_t$. Then, by (2.20) we have

$$
\|\mathbf{y}_0 \widetilde{\mathbf{B}}_{m,0,j}\| \leq C_1 |m^{\lambda_t}| \log^{\nu_t-1} m \leq C_1 V_m. \tag{2.22}
$$

Since the elements of $E(\widetilde{\mathbf{Q}}_i^* \widetilde{\mathbf{Q}}_i)$ are bounded, we have

$$
\begin{aligned}
\left\|\sum_{i=1}^{m} \widetilde{\mathbf{Q}}_i \widetilde{\mathbf{B}}_{m,i,j}\right\| &= \left\{\sum_{i=1}^{m} \widetilde{\mathbf{B}}_{m,i,j}^* E(\widetilde{\mathbf{Q}}_i^* \widetilde{\mathbf{Q}}_i)\widetilde{\mathbf{B}}_{m,i,j}\right\}^{1/2} \\
&\leq C_2 \left\{\sum_{i=1}^{m} (m/i)^{2\operatorname{Re}(\lambda_t)} \log^{2\nu_t-2}(m/i)\right\}^{1/2} \\
&\leq C_2 V_m,
\end{aligned} \tag{2.23}
$$

for all $m$ and some constant $C_2$.

Noticing that $a_{i-1}^{-1}\|\mathbf{y}_{i-1}\|$ is bounded, for $\tau \neq \frac{1}{2}$, we have

$$
\begin{aligned}
\left\|\sum_{i=1}^{m}\left(\frac{i-a_{i-1}}{i}\right)\frac{\mathbf{y}_{i-1}}{a_{i-1}}\mathbf{J}\widetilde{\mathbf{B}}_{m,i,j}\right\| & \\
\leq \sum_{i=1}^{m}\left\|\left(\frac{i-a_{i-1}}{i}\right)\frac{\mathbf{y}_{i-1}}{a_{i-1}}\mathbf{J}\widetilde{\mathbf{B}}_{m,i,j}\right\| & \\
\leq \sum_{i=1}^{m} C_3 i^{-1/2}(|\lambda_t|+1)\left(\frac{m}{i}\right)^{\operatorname{Re}(\lambda_t)}\log^{\nu_t-1}\left(\frac{m}{i}\right) & \\
\leq \frac{C_3}{\sqrt{m}}\sum_{i=1}^{m}\left(\frac{m}{i}\right)^{\tau+1/2}\log^{\nu_t-1}\left(\frac{m}{i}\right) \leq C_3 V_m, &
\end{aligned}
$$

for all $m$ and some constant $C_3$.

Now we estimate this term for the case $\tau = \frac{1}{2}$. We have

$$
\left\|\sum_{i=1}^{m}\left(\frac{i-a_{i-1}}{i}\right)\frac{\mathbf{y}_{i-1}}{a_{i-1}}\mathbf{J}\widetilde{\mathbf{B}}_{m,i,j}\right\|
$$

$$\leq \left\| \sum_{i=1}^{m} \left( \frac{i - a_{i-1}}{i} \right) E \frac{\mathbf{y}_{i-1}}{a_{i-1}} \mathbf{J} \widetilde{\mathbf{B}}_{m,i,j} \right\|$$

$$+ \sum_{i=1}^{m} \left\| \left( \frac{i - a_{i-1}}{i} \right) \left[ \frac{\mathbf{y}_{i-1}}{a_{i-1}} - E \frac{\mathbf{y}_{i-1}}{a_{i-1}} \right] \mathbf{J} \widetilde{\mathbf{B}}_{m,i,j} \right\|.$$

First, we have

$$\left\| \sum_{i=1}^{m} \left( \frac{i - a_{i-1}}{i} \right) E \frac{\mathbf{y}_{i-1}}{a_{i-1}} \mathbf{J} \widetilde{\mathbf{B}}_{m,i,j} \right\|^2$$

$$= \sum_{ik} \left( \frac{E(i - a_{i-1})(k - a_{k-1})}{ik} \right) E \frac{\mathbf{y}_{i-1}}{a_{i-1}} \mathbf{J} \widetilde{\mathbf{B}}_{m,i,j} E \frac{\mathbf{y}_{k-1}}{a_{k-1}} \mathbf{J} \widetilde{\mathbf{B}}_{m,k,j}$$

$$\leq C_3^2 \sum_{i \leq k \leq m} \frac{\|i - a_{i-1}\|^2}{ik} \left( \frac{m}{i} \right)^{1/2} \log^{\nu_t - 1} \left( \frac{m}{i} \right) \left( \frac{m}{k} \right)^{1/2} \log^{\nu_t - 1} \left( \frac{m}{k} \right)$$

$$\leq C_3^2 m \log^{2\nu_t - 2} m \sum_{i \leq k \leq m} \frac{1}{k \sqrt{ik}}$$

$$\leq C_3^2 m \log^{2\nu_t - 1} m.$$

Here we point out a fact that for any $\mu > 1$, there is a constant $C_\mu > 0$ such that

$$E|a_n - n|^\mu \leq C_\mu n^{\mu/2}.$$

This inequality is an easy consequence of the Burkholder inequality. [The Burkholder inequality states that if $X_1, \ldots, X_n$ is a sequence of martingale differences, then for any $p > 1$, there is a constant $C = C(p)$ such that $E|\sum_{i=1}^{m} X_i|^p \leq C_p E(\sum_{i=1}^{n} E(|X_i^2\|\mathcal{F}_{i-1})^{p/2}.]$

By using $\frac{1}{a_{i-1}} = \frac{1}{i} + \frac{i - a_{i-1}}{i a_{i-1}}$ and the above inequality, we have

$$\sum_{i=1}^{m} \left\| \left( \frac{i - a_{i-1}}{i} \right) \left[ \frac{\mathbf{y}_{i-1}}{a_{i-1}} - E \frac{\mathbf{y}_{i-1}}{a_{i-1}} \right] \mathbf{J} \widetilde{\mathbf{B}}_{m,i,j} \right\|$$

$$\leq \sum_{i=1}^{m} \left\| \left( \frac{i - a_{i-1}}{i} \right) \left[ \frac{\mathbf{y}_{i-1} - E\mathbf{y}_{i-1}}{i} \right] \mathbf{J} \widetilde{\mathbf{B}}_{m,i,j} \right\| + \sum_{i=1}^{m} \frac{\|(i - a_{i-1})^2\|}{i^2} \|\mathbf{J} \widetilde{\mathbf{B}}_{m,i,j}\|$$

$$\leq C_3 \sum_{i=1}^{m} \left[ \frac{\|\mathbf{y}_{i-1} - E\mathbf{y}_{i-1}\|}{i^{5/4}} + P(|a_{i-1} - i| \geq i^{3/4}) + \frac{1}{i} \right] \|\mathbf{J} \widetilde{\mathbf{B}}_{m,i,j}\|$$

$$\leq C_3 M_0 \sum_{i=1}^{n_0} \frac{V_i}{i^{5/4}} \left( \frac{m}{i} \right)^{1/2} \log^{\nu_t - 1} \left( \frac{m}{i} \right)$$

$$+ C_3 M \sum_{i=n_0+1}^{m} \frac{V_i}{i^{5/4}} \left(\frac{m}{i}\right)^{1/2} \log^{\nu_t - 1}\left(\frac{m}{i}\right)$$

$$+ C_3 \sum_{i=1}^{m} i^{-1} \left(\frac{m}{i}\right)^{1/2} \log^{\nu_t - 1}\left(\frac{m}{i}\right)$$

$$\leq (C_3 M_0 + C_3 + \varepsilon M)\sqrt{m} \log^{\nu - 1} m.$$

Combining the above four inequalities, we have proved

$$(2.24) \qquad \left\| \sum_{i=1}^{m} \left(\frac{i - a_{i-1}}{i}\right) \frac{\mathbf{y}_{i-1}}{a_{i-1}} \mathbf{J} \widetilde{\mathbf{B}}_{m,i,j} \right\| \leq (C_3 M_0 + C_3 + \varepsilon M) V_m.$$

By (1.3) and the fact that $a_{i-1}^{-1} \|\mathbf{y}_{i-1}\|$ is bounded, we have

$$(2.25) \qquad \begin{aligned}
\sum_{i=1}^{m} & \left\| \left(\frac{i - a_{i-1}}{i}\right) \frac{\mathbf{y}_{i-1}}{a_{i-1}} \mathbf{W}_i \widetilde{\mathbf{B}}_{m,i,j} \right\| \\
& \leq \sum_{i=1}^{m} C_4 i^{-1/2} \alpha_i \left(\frac{m}{i}\right)^{\mathrm{Re}(\lambda_t)} \log^{\nu_t - 1}\left(\frac{m}{i}\right) \\
& \leq C_4 \sqrt{m} \sum_{i=1}^{m} \frac{\alpha_i}{i} \left(\frac{m}{i}\right)^{\tau - 1/2} \log^{\nu_t - 1}\left(\frac{m}{i}\right) \\
& \leq C_4 V_m.
\end{aligned}$$

Next, we show that

$$(2.26) \qquad \begin{aligned}
\sum_{i=1}^{m} & \left\| \frac{\mathbf{y}_{i-1}\mathbf{W}_i - E\mathbf{y}_{i-1}\mathbf{W}_i}{i} \widetilde{\mathbf{B}}_{m,i,j} \right\| \\
& \leq \sum_{i=1}^{m} \left\| \frac{\mathbf{y}_{i-1} - E\mathbf{y}_{i-1}}{i} \mathbf{W}_i \widetilde{\mathbf{B}}_{m,i,j} \right\| \\
& \quad + \sum_{i=1}^{m} \left\| \frac{E(\mathbf{y}_{i-1} - E\mathbf{y}_{i-1})\mathbf{W}_i}{i} \widetilde{\mathbf{B}}_{m,i,j} \right\| \\
& \quad + \sum_{i=1}^{m} \left\| \frac{E\mathbf{y}_{i-1}}{i}(\mathbf{W}_i - E\mathbf{W}_i)\widetilde{\mathbf{B}}_{m,i,j} \right\| \\
& \leq (2\varepsilon M + C_5(2M_0 + 1))V_m.
\end{aligned}$$

By (1.3) and the induction assumption that $\|\mathbf{y}_{i-1} - E\mathbf{y}_{i-1}\| \leq M\sqrt{i}$,

$$\begin{aligned}
\sum_{i=1}^{m} & \left\| \frac{\mathbf{y}_{i-1} - E\mathbf{y}_{i-1}}{i} \mathbf{W}_i \widetilde{\mathbf{B}}_{m,i,j} \right\| \\
& \leq \sum_{i=1}^{n_0} M_0 V_i i^{-1} \alpha_i \left(\frac{m}{i}\right)^{\mathrm{Re}(\lambda_t)} \log^{\nu_t - 1}\left(\frac{m}{i}\right)
\end{aligned}$$

$$+ \sum_{i=n_0+1}^{m} MV_i i^{-1} \alpha_i \left(\frac{m}{i}\right)^{\mathrm{Re}(\lambda_t)} \log^{\nu_t-1}\left(\frac{m}{i}\right)$$

$$\leq (C_5 M_0 + \varepsilon M) V_m.$$

By Jensen's inequality, we have

$$\sum_{i=1}^{m} \left\| E \frac{\mathbf{y}_{i-1} - E\mathbf{y}_{i-1}}{i} \mathbf{W}_i \widetilde{\mathbf{B}}_{m,i,j} \right\| \leq \sum_{i=1}^{m} \left\| \frac{\mathbf{y}_{i-1} - E\mathbf{y}_{i-1}}{i} \mathbf{W}_i \widetilde{\mathbf{B}}_{m,i,j} \right\|$$

$$\leq (C_5 M_0 + \varepsilon M) V_m.$$

The estimate of the third term is given by

$$\sum_{i=1}^{m} \left\| \frac{E\mathbf{y}_{i-1}}{i} (\mathbf{W}_i - E\mathbf{W}_i) \widetilde{\mathbf{B}}_{m,i,j} \right\|$$

$$\leq C_5 \sum_{i=1}^{m} \|\mathbf{W}_i - E\mathbf{W}_i\| \left(\frac{m}{i}\right)^{\mathrm{Re}(\lambda_t)} \log^{\nu_t-1}\left(\frac{m}{i}\right)$$

(2.27)
$$\leq \begin{cases} C_5 \sum_{i=1}^{m} i^{-1/2} \left(\frac{m}{i}\right)^{\mathrm{Re}(\lambda_t)} \log^{\nu_t-1}\left(\frac{m}{i}\right), & \text{if } \tau \neq \frac{1}{2}, \\ C_5 \sum_{i=1}^{m} i^{-1/2} \log^{-1/2}(i+1) \left(\frac{m}{i}\right)^{\mathrm{Re}(\lambda_t)} \log^{\nu_t-1}\left(\frac{m}{i}\right), & \text{if } \tau = \frac{1}{2}, \end{cases}$$

$$\leq C_5 V_m.$$

The above three estimates prove the assertion (2.26).

Substituting (2.22)–(2.26) into (2.21), we obtain

$$\|y_{n,j} - Ey_{n,j}\| \leq (3\varepsilon M + C_1 + C_2 + C_3 + C_4 + C_5 + (C_3 + 2C_5)M_0)V_m \leq MV_m.$$

We complete the proof of (2.15) and thus of (2.12).

Since $\kappa > \tau \vee 1/2$, we may choose $\kappa_1$ such that $\kappa > \kappa_1 > \tau \vee 1/2$. By (2.12), we have

$$\|\mathbf{Y}_n - E\mathbf{Y}_n\|^2 \leq Mn^{2\kappa_1}.$$

From this and the standard procedure of subsequence method, one can show that

$$n^{-\kappa}(\mathbf{Y}_n - E\mathbf{Y}_n) \to 0 \qquad \text{a.s.}$$

To complete the proof of the theorem, it remains to show the replacement of $E\mathbf{Y}_n$ by $n\mathbf{v}$, that is, to show that $\|y_{n,j}\| \leq MV_n$ if (2.13) holds and that $\|y_{n,j}\| = o(n)$ under (1.3). Here the latter is for the convergence with $\kappa = 1$.

Following the lines of the proof for the first conclusion, we need only to change $Ey_{m,j}$ on the left-hand side of (2.21) and replace $E\mathbf{y}_{i-1}\mathbf{W}_i$ on the right-hand side of (2.21) by 0. Checking the proofs of (2.22)–(2.26), we find that the proofs of (2.22)–(2.26) remain true. Therefore, we need only show that

$$
\begin{aligned}
&\sum_{i=1}^{m}\left\|\frac{E\mathbf{y}_{i-1}\mathbf{W}_i}{i}\widetilde{\mathbf{B}}_{m,i,j}\right\| \\
(2.28)\qquad &\leq \sum_{i=1}^{m}\alpha_i\left(\frac{m}{i}\right)^{\mathrm{Re}(\lambda_t)}\log^{\nu_t-1}\left(\frac{m}{i}\right) \\
&=\begin{cases}
\sqrt{m}\displaystyle\sum_{i=1}^{m}\frac{\alpha_i}{\sqrt{i}}\left(\frac{m}{i}\right)^{\tau-1/2}\log^{\nu_t-1}\left(\frac{m}{i}\right)\leq O(V_m), & \text{if (2.13) holds,} \\
m\displaystyle\sum_{i=1}^{m}\frac{\alpha_i}{i}\left(\frac{m}{i}\right)^{\tau-1}\log^{\nu_t-1}\left(\frac{m}{i}\right)\leq\varepsilon m, & \text{if (1.3) holds.}
\end{cases}
\end{aligned}
$$

This completes the proof of this theorem. $\quad\square$

Recall the proof of Theorem 2.2 and note that $\varepsilon$ can be arbitrarily small; with a slight modification to the proof of Theorem 2.2, we have in fact the following corollary.

COROLLARY 2.1. *In addition to the conditions of Theorem 2.2, assume* (2.11) *is true. Then, we have*

$$
(2.29)\qquad \mathbf{y}_{n,-}-E\mathbf{y}_{n,-}=\sum_{i=1}^{n}\widetilde{\mathbf{Q}}_i\widetilde{\mathbf{B}}_{n,i,-}+o_p(V_n),
$$

*where* $\mathbf{y}_{n,-}=(y_{n,2},\ldots,y_{n,K})$ *and* $\widetilde{\mathbf{B}}_{n,i,-}=(\widetilde{\mathbf{B}}_{n,i,2},\ldots,\widetilde{\mathbf{B}}_{n,i,K})$.
*Furthermore, if* (2.13) *is true,* $E\mathbf{y}_{n,-}$ *in* (2.29) *can be replaced by* 0.

PROOF. Checking the proof of Theorem 2.2, one finds that the term estimated in (2.22) is not necessary to appear on the right-hand side of (2.21). Thus, to prove (2.29), it suffices to improve the right-hand sides of (2.24)–(2.26) as to $\varepsilon V_m$. The modification for (2.24) and (2.25) can be done without any further conditions, provided one notices that the vector $\mathbf{y}_{i-1}$ in these inequalities can be replaced by $(0,\mathbf{y}_{i-1,-})$. The details are omitted. To modify (2.26), we first note that (2.27) can be trivially modified to $\varepsilon V_m$ if the condition (2.10) is strengthened to (2.11). The other two estimates for proving (2.26) can be modified easily without any further assumptions. $\quad\square$

Note that

$$
\mathbf{N}_n=\sum_{i=1}^{n}\mathbf{X}_i=\sum_{i=1}^{n}(\mathbf{X}_i-E(\mathbf{X}_i|\mathcal{F}_{i-1}))+\sum_{i=1}^{n}\frac{\mathbf{Y}_{i-1}}{a_{i-1}}.
$$

Since $(\mathbf{X}_i - E(\mathbf{X}_i|\mathcal{F}_{i-1}))$ is a bounded martingale difference sequence, we have

$$n^{-\kappa}\sum_{i=1}^{n}(\mathbf{X}_i - E(\mathbf{X}_i|\mathcal{F}_{i-1})) \to 0 \qquad \text{a.s.}$$

for any $\kappa > 1/2$. Also,

$$n^{-\kappa}\left\|\sum_{i=1}^{n}\mathbf{Y}_{i-1}|i^{-1} - a_{i-1}^{-1}|\right\| \le n^{-\kappa}\sum_{i=1}^{n}i^{-1+\kappa_1}\frac{|a_{i-1}-i|}{i^{\kappa_1}} \to 0 \qquad \text{a.s.}$$

In view of these relations and Theorem 2.2, we have established the following theorem for the strong consistency of $\mathbf{N}_n$.

THEOREM 2.3. *Under the assumptions of Theorem* 2.2, $n^{-\kappa}(\mathbf{N}_n - E\mathbf{N}_n) \to 0$, *a.s. for any* $\kappa > \tau \vee 1/2$. *Also, in the above limit, $E\mathbf{N}_n$ can be replaced by $n\mathbf{v}$ if $\kappa = 1$ or* (2.13) *is true. This implies that $n^{-1}\mathbf{N}_n$ almost surely converges to $\mathbf{v}$, the same limit of $n^{-1}E\mathbf{N}_n$, as $n \to \infty$.*

## 3. Asymptotic normality of $\mathbf{Y}_n$.
In the investigation of the asymptotic normality of the urn composition, we first consider that of $a_n$, the total number of balls in the urn after $n$ stages.

THEOREM 3.1. *Under Assumptions* 2.1–2.3, $n^{-1/2}(a_n - n)$ *is asymptotically normal with mean* 0 *and variance* $\sigma_{11}$, *where* $\sigma_{11} = \sum_{q=1}^{K}\sum_{k=1}^{K}\sum_{l=1}^{K}v_q d_{qkl}$.

PROOF. From Theorems 2.1 and 2.2, we have that $\mathbf{Y}_n/a_n \to \mathbf{v}$ a.s. Similar to (2.8), we have

$$n^{-1}\sum_{i=1}^{n}\operatorname{var}(e_i|\mathcal{F}_{i-1}) \to \sum_{q=1}^{K}\sum_{k=1}^{K}\sum_{l=1}^{K}v_q d_{qkl} \qquad \text{a.s.}$$

Assumption 2.2 implies that $\{e_i - E(e_i|\mathcal{F}_{i-1})\}$ satisfies the Lyapunov condition. From the martingale CLT [see Hall and Heyde (1980)], Assumptions 2.1–2.3 and the fact that

$$a_n - n = 1 + \sum_{i=1}^{n}(e_i - E(e_i|\mathcal{F}_{i-1})),$$

the theorem follows. $\square$

THEOREM 3.2. *Under the assumptions of Theorem* 2.2, $V_n^{-1}(\mathbf{Y}_n - E\mathbf{Y}_n)$ *is asymptotically normal with mean vector* 0 *and variance–covariance matrix* $\mathbf{T}^{-1*}\Sigma\mathbf{T}^{-1}$, *where* $\Sigma$ *is specified later,* $V_n^2 = n$ *if* $\tau < 1/2$ *and* $V_n^2 = n\log^{2\nu-1}n$ *if* $\tau = 1/2$. *Here* $\tau$ *is defined in Assumption* 2.3.
*Also, if* (2.13) *holds, then* $E\mathbf{Y}_n$ *can be replaced by* $n\mathbf{v}$.

PROOF.   To show the asymptotic normality of $\mathbf{Y}_n - E\mathbf{Y}_n$, we only need to show that of $(\mathbf{Y}_n - E\mathbf{Y}_n)\mathbf{T} = \mathbf{y}_n - E\mathbf{y}_n$.

From the proof of Theorem 3.1, we have

$$y_{n,1} - Ey_{n,1} = a_n - n - 1 = \sum_{i=1}^{n}(e_i - E(e_i|\mathcal{F}_{i-1})).$$

From Corollary 2.1, we have

$$\mathbf{y}_{n,-} - E\mathbf{y}_{n,-} = \sum_{i=1}^{n}\widetilde{Q}_i\widetilde{\mathbf{B}}_{n,i,-} + o_p(V_n).$$

Combining the above estimates, we get

$$(3.1)\qquad \mathbf{y}_n - E\mathbf{y}_n = \begin{pmatrix} \displaystyle\sum_{i=1}^{n}(e_i - E(e_2|\mathcal{F}_{i-1})) \\ \displaystyle\sum_{i=1}^{n-1}\widetilde{\mathbf{Q}}_i\widetilde{\mathbf{B}}_{n,i,2} \\ \cdots\cdots\cdots \\ \displaystyle\sum_{i=1}^{n-1}\widetilde{\mathbf{Q}}_i\widetilde{\mathbf{B}}_{n,i,K} \end{pmatrix} + o(V_n).$$

Again, Assumption 2.2 implies the Lyapunov condition. Using the CLT for martingale sequence, as was done in the proof of Theorem 2.3 of Bai and Hu (1999), from (3.1), one can easily show that $V_n^{-1}(\mathbf{y}_n - E\mathbf{y}_n)$ tends to a $K$-variate normal distribution with mean 0 and variance–covariance matrix $\begin{pmatrix} \sigma_{11} & \Sigma_{12} \\ \Sigma_{12}' & \Sigma_{22} \end{pmatrix}$. The variance–covariance matrix $\Sigma_{22}$ of the second to the $K$th elements of $V_n^{-1}(\mathbf{y}_n - E\mathbf{y}_n)$ can be found in (2.17) of Bai and Hu (1999).

By Theorem 3.1, for the case $\tau = 1/2$, $V_n = \sqrt{n}\log^{\nu-1/2}n$, $\sigma_{11} = 0$ and $\Sigma_{12} = 0$. When $\tau < 1/2$, $V_n = \sqrt{n}$, $\sigma_{11} = \sum_{q=1}^{K}\sum_{k=1}^{K}\sum_{l=1}^{K}v_{q}d_{qkl}$. Now, let us find $\Sigma_{12}$. Write $\mathbf{T} = (\mathbf{1}', \mathbf{T}_1, \ldots, \mathbf{T}_s) = (\mathbf{1}', \mathbf{T}_-)$, $\mathbf{T}_j = (\mathbf{t}'_{j1}, \ldots, \mathbf{t}'_{j\nu_j})$ and $\widetilde{\mathbf{B}}_{n,i,-} = \mathbf{T}^{-1}\mathbf{B}_{n,i}\mathbf{T}_- = (\widetilde{\mathbf{B}}_{n,i,2}, \ldots, \widetilde{\mathbf{B}}_{n,i,K})$, where $\mathbf{1} = (1, \ldots, 1)$ throughout this paper. Then the vector $\Sigma_{12}$ is the limit of

$$
\begin{aligned}
n^{-1}&\sum_{i=1}^{n}\operatorname{cov}[(e_i,\mathbf{Q}_i)|\mathcal{F}_{i-1}]\mathbf{T}\widetilde{\mathbf{B}}_{n,i,-}\\
&=n^{-1}\sum_{i=1}^{n}\mathbf{1}\operatorname{cov}[(\mathbf{D}_i'\mathbf{X}_i',\mathbf{X}_i\mathbf{D}_i)|\mathcal{F}_{i-1}]\mathbf{T}\widetilde{\mathbf{B}}_{n,i,-}\\
&=n^{-1}\sum_{i=1}^{n}\mathbf{1}\left(\sum_{q=1}^{K}v_q\mathbf{d}_q+\mathbf{H}^*(\operatorname{diag}(\mathbf{v})-\mathbf{v}^*\mathbf{v})\mathbf{H}\right)\mathbf{T}\widetilde{\mathbf{B}}_{n,i,-}+o_p(1)\\
&=\mathbf{1}\left(\sum_{q=1}^{K}v_q\mathbf{d}_q+\mathbf{H}^*(\operatorname{diag}(\mathbf{v})-\mathbf{v}^*\mathbf{v})\mathbf{H}\right)\mathbf{T}n^{-1}\sum_{i=1}^{n}\widetilde{\mathbf{B}}_{n,i,-}+o_p(1)\\
&=\mathbf{1}\left(\sum_{q=1}^{K}v_q\mathbf{d}_q\right)\mathbf{T}n^{-1}\sum_{i=1}^{n}\widetilde{\mathbf{B}}_{n,i,-}+o_p(1),
\end{aligned}
\tag{3.2}
$$

where the matrices $\mathbf{d}_q$ are defined in (2.4). Here we have used the fact that $\mathbf{1}\mathbf{H}^*(\operatorname{diag}(\mathbf{v})-\mathbf{v}^*\mathbf{v})=\mathbf{1}(\operatorname{diag}(\mathbf{v})-\mathbf{v}^*\mathbf{v})=0$.

By elementary calculation and the definition of $\widetilde{\mathbf{B}}_{n,i,-}$, we get

$$
\begin{aligned}
&n^{-1}\sum_{i=1}^{n}\widetilde{\mathbf{B}}_{n,i,-}^{*}\\
&=\begin{pmatrix}
0 & \cdots & 0\\
n^{-1}\sum_{i=1}^{n}\prod_{j=i+1}^{n}(\mathbf{I}+j^{-1}\mathbf{J}_1) & \cdots & 0\\
\vdots & \cdots & \vdots\\
0 & \cdots & n^{-1}\sum_{i=1}^{n}\prod_{j=i+1}^{n}(\mathbf{I}+j^{-1}\mathbf{J}_s)
\end{pmatrix}.
\end{aligned}
\tag{3.3}
$$

In the $h$th block of the quasi-diagonal matrix

$$
n^{-1}\sum_{i=1}^{n}\prod_{j=i+1}^{n}(\mathbf{I}+j^{-1}\mathbf{J}_h),
$$

the $(g,g+\ell)$-element $(0\le\ell\le\nu_h-1)$ has the approximation

$$
n^{-1}\sum_{i=1}^{n}\frac{1}{\ell!}\left(\frac{n}{i}\right)^{\lambda_h}\log^{\ell}\left(\frac{n}{i}\right)(1+o(1))\to\left(\frac{1}{1-\lambda_h}\right)^{\ell+1}.
\tag{3.4}
$$

Combining (3.2)–(3.4), we get an expression of $\Sigma_{12}$.

Therefore, $n^{-1/2}(\mathbf{y}_n-E\mathbf{y}_n)$ has an asymptotically joint normal distribution with mean 0 and variance–covariance matrix $\Sigma$. Thus, we have shown that

$$
n^{-1/2}(\mathbf{Y}_n-E\mathbf{Y}_n)\to N(0,(\mathbf{T}^{-1})^*\Sigma\mathbf{T}^{-1})
$$

in distribution.

When (2.13) holds, $\mathbf{y}_n - n\mathbf{e}_1$ has the same approximation of the right-hand side of (3.1). Therefore, in the CLT, $E\mathbf{Y}_n$ can be replaced by $n\mathbf{v}$. Then, we complete the proof of the theorem. $\square$

EXAMPLE 3.1. Consider the most important case in application, where $\mathbf{H}$ has a diagonal Jordan form and $\tau < 1/2$. We have

$$\mathbf{T}^{-1}\mathbf{H}\mathbf{T} = \mathbf{J} = \begin{pmatrix} 1 & 0 & \cdots & 0 \\ 0 & \lambda_1 & \cdots & 0 \\ \cdots & \cdots & \cdots & \cdots \\ 0 & 0 & \cdots & \lambda_{K-1} \end{pmatrix},$$

where $\mathbf{T} = (\mathbf{1}', \mathbf{t}_1', \ldots, \mathbf{t}_{K-1}')$. Now let

$$\mathbf{R} = \sum_{j=1}^{K} v_j \mathbf{d}_j + \mathbf{H}^*(\operatorname{diag}(\mathbf{v}) - \mathbf{v}^*\mathbf{v})\mathbf{H}.$$

The variance–covariance matrix $\Sigma = (\sigma_{ij})_{i,j=1}^{K}$ has the following simple form:
$\sigma_{11} = \mathbf{1}\mathbf{R}\mathbf{1}' = \sum_{q=1}^{K}\sum_{k=1}^{K}\sum_{l=1}^{K} v_q d_{qkl}$, $\sigma_{1j} = (1-\lambda_{j-1})^{-1}\mathbf{1}\mathbf{R}\mathbf{t}_{j-1}' = (1-\lambda_{j-1})^{-1}\sum_{k=1}^{K} v_k \mathbf{1}\mathbf{d}_k\mathbf{t}_{j-1}'$, $j = 2, \ldots, K$, and

$$\sigma_{ij} = (1 - \lambda_{i-1} - \bar{\lambda}_{j-1})^{-1}(\mathbf{t}_{i-1}^*)'\mathbf{R}\mathbf{t}_{j-1}'.$$

**4. Asymptotic normality of $\mathbf{N}_n$.** Now, $\mathbf{N}_n = (N_{n1}, \ldots, N_{nK})$, where $N_{nk}$ is the number of times a type-$k$ ball is drawn in the first $n$ draws:

$$\mathbf{N}_n = (N_{n1}, \ldots, N_{nK}) = \mathbf{N}_{n-1} + \mathbf{X}_n = \sum_{i=1}^{n} \mathbf{X}_i,$$

where the vectors $\mathbf{X}_i$ are defined as follows: If a type-$k$ ball is drawn in the $i$th stage, then define the draw outcome $\mathbf{X}_i$ as the vector whose $k$th component is 1 and all others are 0. Therefore $\mathbf{1}\mathbf{X}_i' = 1$ and $\mathbf{1}\mathbf{N}_n' = n$. We shall consider the limiting property of $\mathbf{N}_n$.

THEOREM 4.1 (for the EPU urn). *Under the assumptions of Corollary* 2.1, $V_n^{-1}(\mathbf{N}_n - E\mathbf{N}_n)$ *is asymptotically normal with mean vector 0 and variance–covariance matrix* $\mathbf{T}^{-1*}\widetilde{\Sigma}\mathbf{T}^{-1}$, *where* $\widetilde{\Sigma}$ *is specified later,* $V_n^2 = n$ *if* $\tau < 1/2$ *and* $V_n^2 = n\log^{2\nu-1} n$ *if* $\tau = 1/2$. *Here* $\tau$ *is defined in Assumption* 2.3.

*Furthermore, if* (2.13) *holds, then* $E\mathbf{N}_n$ *can be replaced by* $n\mathbf{v}$.

PROOF. At first we have

$$(4.1) \qquad \begin{aligned} \mathbf{N}_n &= \sum_{i=1}^{n} (\mathbf{X}_i - E(\mathbf{X}_i|\mathcal{F}_{i-1})) + \sum_{i=1}^{n} E(\mathbf{X}_i|\mathcal{F}_{i-1}) \\ &= \sum_{i=1}^{n} (\mathbf{X}_i - \mathbf{Y}_{i-1}/a_{i-1}) + \sum_{i=0}^{n-1} \mathbf{Y}_i/a_i. \end{aligned}$$

For simplicity, we consider the asymptotic distribution of $\mathbf{N}_n\mathbf{T}$. Since the first component of $\mathbf{N}_n\mathbf{T}$ is a nonrandom constant $n$, we only need consider the other $K-1$ components. From (2.29) and (4.1), we get

$$
\begin{aligned}
(4.2)\quad \mathbf{N}_n\mathbf{T}_- &= \sum_{i=1}^n (\mathbf{X}_i - \mathbf{Y}_{i-1}/a_{i-1})\mathbf{T}_- + \sum_{i=0}^{n-1} \mathbf{y}_{i,-}/a_i \\
&= \sum_{i=1}^n (\mathbf{X}_i - \mathbf{Y}_{i-1}/a_{i-1})\mathbf{T}_- \\
&\quad + \sum_{i=0}^{n-1} \mathbf{y}_{i,-}/(i+1) + \sum_{i=0}^{n-1} \frac{\mathbf{y}_{i,-}}{a_i}\left(\frac{i+1-a_i}{i+1}\right) \\
&= \sum_{i=1}^n (\mathbf{X}_i - \mathbf{Y}_{i-1}/a_{i-1})\mathbf{T}_- + \sum_{i=0}^{n-1} \frac{\mathbf{y}_{i,-}}{a_i}\left(\frac{i+1-a_i}{i+1}\right) + \mathbf{y}_{0,-} \\
&\quad + \sum_{i=1}^{n-1} \frac{1}{i+1}\left[\sum_{j=1}^i \widetilde{\mathbf{Q}}_j\widetilde{\mathbf{B}}_{i,j,-} + E\mathbf{y}_{i,-} + o_p(V_i)\right] \\
&= \sum_{i=1}^n (\mathbf{X}_i - \mathbf{Y}_{i-1}/a_{i-1})\mathbf{T}_- \\
&\quad + \sum_{j=1}^{n-1} \widetilde{\mathbf{Q}}_j\left(\sum_{i=j}^{n-1}\frac{1}{i+1}\widetilde{\mathbf{B}}_{i,j,-}\right) + \sum_{i=1}^{n-1}\frac{E\mathbf{y}_{i,-}}{i+1} + o(V_n) \\
&= \sum_{i=1}^n (\mathbf{X}_i - \mathbf{Y}_{i-1}/a_{i-1})\mathbf{T} + \sum_{j=1}^{n-1} \widetilde{\mathbf{Q}}_j\widehat{\mathbf{B}}_{n,j,-} + \sum_{i=1}^{n-1}\frac{E\mathbf{y}_{i,-}}{i+1} + o_p(V_n),
\end{aligned}
$$

where $\widetilde{\mathbf{B}}_{i,j} = \mathbf{T}^{-1}\mathbf{G}_{j+1}\cdots\mathbf{G}_i\mathbf{T}$, $\widehat{\mathbf{B}}_{n,j} = \sum_{i=j}^{n-1}\frac{1}{i+1}\widetilde{\mathbf{B}}_{i,j}$, the matrices with a minus sign in subscript denote the submatrices of the last $K-1$ columns of their corresponding mother matrices. Here, in the fourth equality, we have used the fact that $\sum_{i=0}^{n-1}\frac{\mathbf{y}_{i,-}}{a_i}(\frac{i+1-a_i}{i+1}) = o_p(\sqrt{n})$ which can be proven by the same approach as showing (2.24) and (2.28).

In view of (4.2), we only have to consider the asymptotic distribution of the martingale

$$
\mathbf{U} = \sum_{i=1}^n (\mathbf{X}_i - \mathbf{Y}_{i-1}/a_{i-1})\mathbf{T}_- + \sum_{j=1}^{n-1} \widetilde{\mathbf{Q}}_j\widehat{\mathbf{B}}_{n,j,-}.
$$

We now estimate the asymptotic variance–covariance matrix of $V_n^{-1}\mathbf{U}$. For this end, we need only consider the limit of

$$
\begin{aligned}
(4.3)\quad \widetilde{\Sigma}_n = V_n^{-2}\bigg[&\sum_{j=1}^n E(\mathbf{q}_j^*\mathbf{q}_j|\mathcal{F}_{j-1}) + \sum_{j=1}^{n-1} E(\mathbf{q}_j^*\widetilde{\mathbf{Q}}_j\widehat{\mathbf{B}}_{n,j,-}|\mathcal{F}_{j-1}) \\
&+ \sum_{j=1}^{n-1} E(\widehat{\mathbf{B}}_{n,j,-}^*\widetilde{\mathbf{Q}}_j^*\mathbf{q}_j|\mathcal{F}_{j-1}) + \sum_{j=1}^{n-1} E(\widehat{\mathbf{B}}_{n,j,-}^*\widetilde{\mathbf{R}}_j\widehat{\mathbf{B}}_{n,j,-}|\mathcal{F}_{j-1})\bigg],
\end{aligned}
$$

where $\mathbf{q}_j = (\mathbf{X}_j - \mathbf{Y}_{j-1}/a_{j-1})\mathbf{T}_-$ and $\widetilde{\mathbf{R}}_j = E(\widetilde{\mathbf{Q}}_j^*\widetilde{\mathbf{Q}}_j|\mathcal{F}_{j-1}) = \mathbf{T}^*\mathbf{R}_j\mathbf{T}$.

From Theorem 3.1, we know that

$$E(\mathbf{q}_j^*\mathbf{q}_j|\mathcal{F}_{j-1}) \to \mathbf{T}^*_-(\operatorname{diag}(\mathbf{v}) - \mathbf{v}^*\mathbf{v})\mathbf{T}_- = \mathbf{T}^*_-\operatorname{diag}(\mathbf{v})\mathbf{T}_- \qquad \text{as } j \to \infty,$$

since $\mathbf{v}\mathbf{T}_- = 0$. This estimate implies that

$$
\begin{aligned}
(4.4) \quad & V_n^{-2}\sum_{j=1}^n E(\mathbf{q}_j^*\mathbf{q}_j|\mathcal{F}_{j-1}) \\
& \to \widetilde{\Sigma}_1 = \begin{cases} \mathbf{T}^*_-\operatorname{diag}(\mathbf{v})\mathbf{T}_-, & \text{if } \tau < 1/2, \\ 0, & \text{if } \tau = 1/2, \end{cases} \quad \text{as } j \to \infty.
\end{aligned}
$$

Because $\widetilde{\mathbf{Q}}_j = [\mathbf{X}_j\mathbf{D}_j - (\mathbf{Y}_{j-1}/a_{j-1})\mathbf{H}_j]\mathbf{T}$,

$$
\begin{aligned}
(4.5) \quad & V_n^{-2}\sum_{j=1}^{n-1} E(\mathbf{q}_j^*\widetilde{\mathbf{Q}}_j\widehat{\mathbf{B}}_{n,j,-}|\mathcal{F}_{j-1}) \\
& = V_n^{-2}\sum_{j=1}^{n-1}\mathbf{T}^*_- E(\mathbf{X}_j - \mathbf{Y}_{j-1}/a_{j-1})^* \\
& \qquad \times [\mathbf{X}_j\mathbf{D}_j - (\mathbf{Y}_{j-1}/a_{j-1})\mathbf{H}_j|\mathcal{F}_{j-1}]\mathbf{T}\widehat{\mathbf{B}}_{n,j,-} \\
& = \mathbf{T}^*_-\operatorname{diag}(\mathbf{v})\mathbf{H}\mathbf{T}\left(V_n^{-2}\sum_{j=1}^{n-1}\widehat{\mathbf{B}}_{n,j,-}\right) + o(1).
\end{aligned}
$$

From (2.18), we have

$$(4.6) \quad n^{-1}\sum_{j=1}^{n-1}\widehat{\mathbf{B}}_{n,j} = n^{-1}\sum_{j=1}^{n-1}\sum_{i=j}^{n-1}\frac{1}{i+1}(\mathbf{I} + (j+1)^{-1}\mathbf{J})\cdots(\mathbf{I} + i^{-1}\mathbf{J}).$$

Based on (2.18)–(2.20), we have that the $(h, h+\ell)$-element of the block matrix

$$n^{-1}\sum_{j=1}^n\sum_{i=j}^{n-1}\frac{1}{i+1}(\mathbf{I} + (j+1)^{-1}\mathbf{J}_t)\cdots(\mathbf{I} + i^{-1}\mathbf{J}_t)$$

has a limit obtained by

$$
\begin{aligned}
(4.7) \quad & n^{-1}\sum_{j=1}^n\sum_{i=j}^{n-1}\frac{1}{i+1}\frac{1}{\ell!}\left(\frac{j}{i}\right)^{-\lambda_t}\log^\ell\left(\frac{i}{j}\right)(1+o(1)) \\
& \to \frac{1}{\ell!}\int_0^1\int_0^u u^{-1}\left(\frac{v}{u}\right)^{-\lambda_t}\log^\ell\left(\frac{u}{v}\right)du\,dv \\
& = \left(\frac{1}{1-\lambda_t}\right)^{\ell+1}.
\end{aligned}
$$

Substituting this into (4.6) and then (4.5), when $V_n^2 = n$, we obtain that

$$V_n^{-2}\sum_{j=1}^{n-1} E(\mathbf{q}_j^*\widetilde{\mathbf{Q}}_j\widehat{\mathbf{B}}_{n,j,-}|\mathcal{F}_{j-1}) \to \widetilde{\Sigma}_2 = \mathbf{T}^*_-\operatorname{diag}(\mathbf{v})\mathbf{H}\mathbf{T}\widetilde{\mathbf{J}},$$

where $\widetilde{\mathbf{J}}$ is a $K \times (K-1)$ matrix whose first row is 0 and the rest is a block diagonal matrix, the $t$-block is $\nu_t \times \nu_t$ and its $(h, h+\ell)$-element is given by the right-hand side of (4.7). The matrix $\widetilde{\Sigma}_2$ is obviously 0 when $V_n^2 = n \log^{2\nu-1} n$.

Note that the third term in (4.3) is the complex conjugate transpose of the second term; thus we have also got the limit of the third term, that is, $\widetilde{\Sigma}_2^*$.

Now, we compute the limit $\widetilde{\Sigma}_3$ of the fourth term in (4.3). By Assumption 2.2, the matrices $\mathbf{R}_j$ in (4.3) converge to $\mathbf{R}$. Then, the fourth term in (4.3) can be approximated by

$$(4.8) \qquad \left[ V_n^{-2} \sum_{j=1}^{n} \sum_{i=j+1}^{n-1} \frac{1}{i} \prod_{r=j+1}^{i} (\mathbf{I} + r^{-1}\mathbf{J}_g^*) \mathbf{T}_g^* \mathbf{R} \mathbf{T}_h \right.$$
$$\left. \times \sum_{i=j+1}^{n-1} \frac{1}{i} \prod_{r=j+1}^{n} (\mathbf{I} + r^{-1}\mathbf{J}_h) \right]_{g,h=1}^{s}.$$

Similar to (4.7), we can show that the $(w, t)$-element of the $(g, h)$-block of the matrix in (4.8) is approximately

$$(4.9) \qquad V_n^{-2} \sum_{w'=0}^{w-1} \sum_{t'=0}^{t-1} \sum_{j=1}^{n-1} \sum_{i=j}^{n-1} \sum_{m=j}^{n-1} \frac{(i/j)^{\bar{\lambda}_g} (m/j)^{\lambda_h} \log^{w'}(i/j) \log^{t'}(m/j)}{(i+1)(m+1)(w')!(t')!}$$
$$\times [\mathbf{T}_g^* \mathbf{R} \mathbf{T}_h]_{(w-w', t-t')},$$

where $[\mathbf{T}_g^* \mathbf{R} \mathbf{T}_h]_{(w', t')}$ is the $(w', t')$-element of the matrix $[\mathbf{T}_g^* \mathbf{R} \mathbf{T}_h]$. Here, strictly speaking, in the numerator of (4.9), there should be factors $\psi(i, j, w')$ and $\psi(m, j, t')$. Since for any $j_0$, the total contributions of terms with $j \le j_0$ is $o(1)$ and the $\psi$'s tend to 1 as $j \to \infty$, we may replace the $\psi$'s by 1.

For fixed $w, w', t$ and $t'$, if $\lambda_g \ne \lambda_h$ or $\mathrm{Re}(\lambda_g) < 1/2$, we have

$$(4.10) \qquad \frac{1}{n} \sum_{j=1}^{n-1} \sum_{i=j}^{n-1} \sum_{m=j}^{n-1} \frac{(i/j)^{\bar{\lambda}_g} (m/j)^{\lambda_h} \log^{w'}(i/j) \log^{t'}(m/j)}{(i+1)(m+1)(w')!(t')!}$$
$$\to \sum_{\ell=0}^{w'} \frac{(t'+\ell)!}{\ell!(t')!} (1-\bar{\lambda}_g)^{-(w'-\ell+1)} (1-\bar{\lambda}_g-\lambda_h)^{-(t'+\ell+1)}$$
$$+ \sum_{\ell=0}^{t'} \frac{(w'+\ell)!}{\ell!(w')!} (1-\lambda_h)^{-(t'-\ell+1)} (1-\bar{\lambda}_g-\lambda_h)^{-(w'+\ell+1)}.$$

Thus, when $\tau < 1/2$, if we split $\widetilde{\Sigma}_3$ into blocks, then the $(w, t)$-element of the $(g, h)$-block $\Sigma_{g,h}$ $(\nu_g \times \nu_h)$ of $\widetilde{\Sigma}_3$ is given by

$$
\begin{aligned}
(4.11) \quad & \sum_{w'=0}^{w-1} \sum_{t'=0}^{t-1} \Bigg[ \sum_{\ell=0}^{w'} \frac{(t'+\ell)!}{\ell!(t')!} (1 - \bar{\lambda}_g)^{-(w'-\ell+1)} (1 - \bar{\lambda}_g - \lambda_h)^{-(t'+\ell+1)} \\
& \qquad\qquad + \sum_{\ell=0}^{t'} \frac{(w'+\ell)!}{\ell!(w')!} (1 - \lambda_h)^{-(t'-\ell+1)} (1 - \bar{\lambda}_g - \lambda_h)^{-(w'+\ell+1)} \Bigg] \\
& \qquad\qquad \times [\mathbf{T}_g^* \mathbf{R} \mathbf{T}_h]_{(w-w', t-t')}.
\end{aligned}
$$

When $\tau = 1/2$, $\Sigma_{g,h} = 0$ if $\lambda_g \neq \lambda_h$ or if $\mathrm{Re}(\lambda_g) < 1/2$. Now, we consider $\Sigma_{g,h}$ with $\lambda_g = \lambda_h$ and $\mathrm{Re}(\lambda_g) = 1/2$. If $w' + t' < 2\nu - 2$, then

$$
\begin{aligned}
& \sum_{j=1}^{n-1} \sum_{i=j}^{n-1} \sum_{\ell=j}^{n-1} \frac{(i/j)^{\bar{\lambda}_g} (\ell/j)^{\lambda_g} \log^{w'}(i/j) \log^{t'}(\ell/j)}{(i+1)(\ell+1)(w')!(t')!} \\
& \quad \le \sum_{j=1}^{n-1} \sum_{i=j}^{n-1} \sum_{\ell=j}^{n-1} \frac{\log^{w'+t'} n}{j\sqrt{i\ell}(w')!(t')!} \le n \log^{w'+t'+1} n = o(V_n^2).
\end{aligned}
$$

When $w' = t' = \nu - 1$ which implies $w = t = \nu = \nu_g = \nu_h$, by Abelian summation, we have

$$
\begin{aligned}
(4.12) \quad & V_n^{-2} \sum_{j=1}^{n-1} \sum_{i=j}^{n-1} \sum_{\ell=j}^{n-1} \frac{(i/j)^{\bar{\lambda}_g} (\ell/j)^{\lambda_g} \log^{\nu-1}(i/j) \log^{\nu-1}(\ell/j)}{(i+1)(\ell+1)[(\nu-1)!]^2} \\
& \quad \to |\lambda_g|^{-2} [(\nu-1)!]^{-2} (2\nu-1)^{-1}.
\end{aligned}
$$

Hence, for this case, $\Sigma_{g,h}$ has only one nonzero element which is the one on the right-lower corner of $\Sigma_{g,h}$ and given by

$$
(4.13) \quad |\lambda_g|^{-2} [(\nu-1)!]^{-2} (2\nu-1)^{-1} [\mathbf{T}_g^* \mathbf{R} \mathbf{T}_h]_{(1,1)}.
$$

Combining (4.3), (4.4), (4.7), (4.11) and (4.12), we obtain an expression of $\widetilde{\Sigma}$. $\qquad\square$

Now we consider one of the most important special cases, where the matrix $\mathbf{H}$ has a diagonal Jordan form and $\tau < 1/2$.

COROLLARY 4.1. *Suppose the assumptions of Corollary* 2.1 *hold with* $\tau < 1/2$ *and*

$$
\mathbf{T}^{-1} \mathbf{H} \mathbf{T} = \mathbf{J} = \begin{pmatrix} 1 & 0 & \cdots & 0 \\ 0 & \lambda_1 & \cdots & 0 \\ \cdots & \cdots & \cdots & \cdots \\ 0 & 0 & \cdots & \lambda_{K-1} \end{pmatrix},
$$

*where* $\mathbf{T} = (\mathbf{1}', \mathbf{t}_1', \dots, \mathbf{t}_{K-1}')$. *Now let*

$$a_{ij} = (\mathbf{t}_{i-1}^*)'(\operatorname{diag}(\mathbf{v}) - \mathbf{v}^*\mathbf{v})\mathbf{t}_{j-1}',$$

$$b_{ij} = \lambda_{j-1}(1 - \lambda_{j-1})^{-1}(\mathbf{t}_{i-1}^*)'(\operatorname{diag}(\mathbf{v}) - \mathbf{v}^*\mathbf{v})\mathbf{t}_{j-1}'$$

*and*

$$c_{ij} = [(1 - \bar{\lambda}_{i-1})^{-1} + (1 - \lambda_{j-1})^{-1}](1 - \bar{\lambda}_{i-1} - \lambda_{j-1})^{-1}(\mathbf{t}_{i-1}^*)'\mathbf{R}\mathbf{t}_{j-1}',$$

*for* $i, j = 2, \dots, K$. *Then* $n^{-1/2}(\mathbf{N}_n - E\mathbf{N}_n)$ *is asymptotically normal with mean vector* $0$ *and variance–covariance matrix* $(\mathbf{T}^{-1})^*\widetilde{\Sigma}\mathbf{T}^{-1}$, *where* $\widetilde{\Sigma} = (\widetilde{\sigma}_{ij})_{i,j=1}^K$ *has the following simple form:*

$$\widetilde{\sigma}_{11} = \widetilde{\sigma}_{1j} = \widetilde{\sigma}_{i1} = 0 \quad \text{and} \quad \widetilde{\sigma}_{ij} = a_{ij} + b_{ij} + \bar{b}_{ji} + c_{ij}$$

*for* $i, j = 2, \dots, K$.

## 5. Applications.

5.1. *Adaptive allocation rules associated with covariates.* In clinical trials, it is usual that the probability of success (here we assume that the subject response is dichotomous) may depend upon some observable covariates on the patients, that is, $p_{ik} = p_k(\xi_i)$, where $\xi_i$ are covariates observed on the patient $i$ and the result of the treatment at the $i$th stage. Here $p_{ik} = P(T_i = 1 | X_i = k, \xi_i)$, for $i = 1, \dots, n$ and $k = 1, \dots, K$, where $X_i = k$ indicates that a type-$k$ ball is drawn at the $i$th stage and $T_i = 1$ if the response of the subject $i$ is a success and 0 otherwise. Thus, for a given $\xi_i$, the addition rule could be $\mathbf{D}(\xi_i)$ and the generating matrices $\mathbf{H}_i = \mathbf{H}(\xi_i) = E\mathbf{D}(\xi_i)$.

Assume that $\xi_1, \dots, \xi_n$ are i.i.d. random vectors and let $\mathbf{H} = E\mathbf{H}(\xi_1)$. The asymptotic properties of the urn composition $\mathbf{Y}_n$ are considered by Bai and Hu (1999). Based on the results in Sections 2 and 4, we can get the corresponding asymptotic results of the allocation number of patients $\mathbf{N}_n$. Here we illustrate the results by considering the case $K = 2$.

Consider the generalized play-the-winner rule [Bai and Hu (1999)] and let $E(p_k(\xi_i)) = p_k$, $k = 1, 2$. Then the addition rule matrices are denoted by

$$\mathbf{D}(\xi_i) = \begin{pmatrix} d_1(\xi_i) & 1 - d_1(\xi_i) \\ 1 - d_2(\xi_i) & d_2(\xi_i) \end{pmatrix} \quad \text{and} \quad \mathbf{H} = \begin{pmatrix} p_1 & q_1 \\ q_2 & p_2 \end{pmatrix},$$

where $0 \le d_k(\xi_i) \le 1$ and $q_k = 1 - p_k$ for $k = 1, 2$.

It is easy to see that $\lambda = 1$, $\lambda_1 = p_1 + p_2 - 1$, $\tau = \max(0, \lambda_1)$ and $\mathbf{v} = (q_2/(q_1 + q_2), q_1/(q_1 + q_2))$. Further, we have

$$\mathbf{R} = \frac{(a_1 q_2 + a_2 q_1)(q_1 + q_2) + q_1 q_2 (p_1 - q_2)^2}{(q_1 + q_2)^2} \begin{pmatrix} 1 & -1 \\ -1 & 1 \end{pmatrix},$$

$$\mathbf{T} = \begin{pmatrix} 1 & q_1 \\ 1 & -q_2 \end{pmatrix} \quad \text{and} \quad \mathbf{T}^{-1} = \frac{1}{q_1 + q_2} \begin{pmatrix} q_2 & q_1 \\ 1 & -1 \end{pmatrix},$$

where $a_k = \text{Var}(d_k(\xi_1))$. For the case $\tau < 1/2$, we have that $V_n = n$ and the values corresponding to Corollary 4.1 are

$$a_{22} = q_1 q_2, \qquad b_{22} = \frac{(1 - q_1 - q_2)q_1 q_2}{q_1 + q_2},$$

$$c_{22} = \frac{2[(a_1 q_2 + a_2 q_1)(q_1 + q_2) + q_1 q_2 (p_1 - q_2)^2]}{(q_1 + q_2)(1 - 2(p_1 + p_2 - 1))}.$$

So

$$\widetilde{\sigma}_{22} = q_1 q_2 + \frac{2(1 - q_1 - q_2)q_1 q_2}{q_1 + q_2} + \frac{2[(a_1 q_2 + a_2 q_1)(q_1 + q_2) + q_1 q_2 (p_1 - q_2)^2]}{(q_1 + q_2)(1 - 2(p_1 + p_2 - 1))}.$$

From Theorem 2.3 and Corollary 4.1, we have

$$n^\delta \left( \frac{\mathbf{N}_n}{n} - \mathbf{v} \right) \to 0 \qquad \text{a.s. for any } \delta < 1/2 \quad \text{and} \quad n^{1/2} \left( \frac{\mathbf{N}_n}{n} - \mathbf{v} \right) \to N(0, \Sigma_1)$$

in distribution, where

$$\Sigma_1 = (\mathbf{T}^{-1})^* \begin{pmatrix} 0 & 0 \\ 0 & \widetilde{\sigma}_{22} \end{pmatrix} \mathbf{T}^{-1} = \frac{\widetilde{\sigma}_{22}}{(q_1 + q_2)^2} \begin{pmatrix} 1 & -1 \\ -1 & 1 \end{pmatrix}.$$

For the randomized play-the-winner rule [Wei and Durham (1978)], we have $a_k = p_k q_k$, $k = 1, 2$. Then we have

$$\widetilde{\sigma}_{22} = \frac{(5 - 2(q_1 + q_2))q_1 q_2}{2(q_1 + q_2) - 1}.$$

This result agrees with that of Matthews and Rosenberger (1997).

For the case $\tau = 1/2$, $V_n = n \log n$ and the value corresponding to (4.11) is

$$\widetilde{\sigma}_{22} = 4[(a_1 q_2 + a_2 q_1)(q_1 + q_2) + q_1 q_2 (p_1 - q_2)^2].$$

We have

$$(n \log n)^{-1/2} (\mathbf{N}_n - n\mathbf{v}) \to N(0, \Sigma_2)$$

in distribution, where

$$\Sigma_2 = \frac{4[(a_1 q_2 + a_2 q_1)(q_1 + q_2) + q_1 q_2 (p_1 - q_2)^2]}{(q_1 + q_2)^2} \begin{pmatrix} 1 & -1 \\ -1 & 1 \end{pmatrix}.$$

For the case of the randomized play-the-winner rule, we have

$$\Sigma_2 = \frac{4 q_1 q_2}{(q_1 + q_2)^2} \begin{pmatrix} 1 & -1 \\ -1 & 1 \end{pmatrix}.$$

5.2. *Clinical trials with time trend in adaptive designs.* Time trends are present in many sequential experiments. Hu and Rosenberger (2000) have studied time trend in adaptive designs and applied to a neurophysiology experiment. It is important to know the asymptotic behavior of the allocation number of patients in these cases.

In Section 5.1, $p_{ik} = P(T_i = 1 | X_i = k)$, where $X_i = k$ if the $k$th element of $\mathbf{X}_i$ is 1. There may be a drift in patient characteristics over time, for example, $\lim_{i \to \infty} p_{ik} = p_k$ [Hu and Rosenberger (2000)]. Then the results in Sections 2, 3 and 4 are applicable here. For the case $K = 2$, we can get similar results as in Section 5.1.

The results in this paper may also apply for GFU model with homogeneous generating matrix with a general Jordan form as well as $\tau = 1/2$. In these cases, the results of Smythe (1996) are not applicable.

5.3. *Urn models for multi-arm clinical trials.* For multi-arm clinical trials, Wei (1979) proposed the following urn model (as an extension of the randomized play-the-winner rule of two treatments): Starting from $\mathbf{Y}_0 = (Y_{01}, \ldots, Y_{0K})$, when a type $k$ splits (randomly from the urn), we assign the patient to the treatment $k$ and observe the patient's response. A success on treatment $k$ adds a ball of type $k$ to the urn and a failure on treatment $k$ adds $1/(K-1)$ ball for each of the other $K-1$ types. Let $p_k$ be the probability of success of treatment $k$, $k = 1, 2, \ldots, K$, and $q_k = 1 - p_k$. The generating matrix for this urn model is

$$\mathbf{H} = \begin{pmatrix} p_1 & (K-1)^{-1}q_1 & \cdots & (K-1)^{-1}q_1 \\ (K-1)^{-1}q_2 & p_2 & \cdots & (K-1)^{-1}q_2 \\ \cdots & \cdots & \cdots & \cdots \\ (K-1)^{-1}q_K & (K-1)^{-1}q_K & \cdots & p_K \end{pmatrix}.$$

The asymptotic properties of $\mathbf{Y}_n$ can be obtained from Athreya and Karlin (1968) and Bai and Hu (1999). From Theorem 4.1 in Section 4, we obtain the asymptotic normality of $\mathbf{N}_n$ and its asymptotic variance.

Recently, Bai, Hu and Shen (2002) proposed an urn model which adds balls depending on the success probabilities of each treatment. Write $\mathbf{N}_n = (N_{n1}, \ldots, N_{nK})$ and $\mathbf{S}_n = (S_{n1}, \ldots, S_{nK})$, where $N_{nk}$ denotes the number of times that the $k$th treatment is selected in the first $n$ stages, and $S_{nk}$ denotes the number of successes of the $k$th treatment in the $N_{nk}$ trials, $k = 1, \ldots, K$. Define $\mathbf{R}_n = (R_{n1}, \ldots, R_{nK})$ and $M_n = \sum_{k=1}^{K} R_{nk}$, where $R_{n,k} = \frac{S_{nk}+1}{N_{nk}+1}$, $k = 1, \ldots, K$. The generating matrices are

$$\mathbf{H}_{i+1} = \begin{pmatrix} p_1 & \dfrac{R_{i2}}{M_i - R_{i1}}q_1 & \cdots & \dfrac{R_{iK}}{M_i - R_{i1}}q_1 \\ \dfrac{R_{i1}}{M_i - R_{i2}}q_2 & p_2 & \cdots & \dfrac{R_{iK}}{M_i - R_{i2}}q_2 \\ \cdots & \cdots & \cdots & \cdots \\ \dfrac{R_{i1}}{M_i - R_{iK}}q_K & \dfrac{R_{i2}}{M_i - R_{iK}}q_K & \cdots & p_K \end{pmatrix}.$$

In this case, $\mathbf{H}_i$ are random matrices and converge to

$$\mathbf{H} = \begin{pmatrix} p_1 & \dfrac{p_2}{M-p_1}q_1 & \cdots & \dfrac{p_K}{M-p_1}q_1 \\ \dfrac{p_1}{M-p_2}q_2 & p_2 & \cdots & \dfrac{p_K}{M-p_2}q_2 \\ \cdots & \cdots & \cdots & \cdots \\ \dfrac{p_1}{M-p_K}q_K & \dfrac{p_2}{M-p_K}q_K & \cdots & p_K \end{pmatrix}$$

almost surely, where $M = p_1 + \cdots + p_K$.

Bai, Hu and Shen (2002) considered the convergences of $\mathbf{Y}_n/n$ and $\mathbf{N}_n/n$. The asymptotic distributions of $\mathbf{Y}_n$ and $\mathbf{N}_n$ can be obtained from Theorems 3.2 and 4.1 in this paper. From Lemma 3 of Bai, Hu and Shen (2002) we have $\alpha_i = o(i^{-1/4})$ almost surely, so the condition (1.3) is satisfied.

**Acknowledgments.** Special thanks go to anonymous referees for the constructive comments, which led to a much improved version of the paper. We would also like to thank Professor W. F. Rosenberger for his valuable discussions which led to the problem of this paper.

## REFERENCES

ALTMAN, D. G. and ROYSTON, J. P. (1988). The hidden effect of time. *Statist. Med.* **7** 629–637.

ANDERSEN, J., FARIES, D. and TAMURA, R. N. (1994). Randomized play-the-winner design for multi-arm clinical trials. *Comm. Statist. Theory Methods* **23** 309–323.

ATHREYA, K. B. and KARLIN, S. (1967). Limit theorems for the split times of branching processes. *Journal of Mathematics and Mechanics* **17** 257–277. MR216592

ATHREYA, K. B. and KARLIN, S. (1968). Embedding of urn schemes into continuous time branching processes and related limit theorems. *Ann. Math. Statist.* **39** 1801–1817. MR232455

BAI, Z. D. and HU, F. (1999). Asymptotic theorem for urn models with nonhomogeneous generating matrices. *Stochastic Process. Appl.* **80** 87–101. MR1670107

BAI, Z. D., HU, F. and SHEN, L. (2002). An adaptive design for multi-arm clinical trials. *J. Multivariate Anal.* **81** 1–18. MR1901202

COAD, D. S. (1991). Sequential tests for an unstable response variable. *Biometrika* **78** 113–121. MR1118236

FLOURNOY, N. and ROSENBERGER, W. F., eds. (1995). *Adaptive Designs.* IMS, Hayward, CA. MR1477667

FREEDMAN, D. (1965). Bernard Friedman's urn. *Ann. Math. Statist.* **36** 956–970. MR177432

GOUET, R. (1993). Martingale functional central limit theorems for a generalized Pólya urn. *Ann. Probab.* **21** 1624–1639. MR1235432

HALL, P. and HEYDE, C. C. (1980). *Martingale Limit Theory and Its Application.* Academic Press, London. MR624435

HOLST, L. (1979). A unified approach to limit theorems for urn models. *J. Appl. Probab.* **16** 154–162. MR520945

HU, F. and ROSENBERGER, W. F. (2000). Analysis of time trends in adaptive designs with application to a neurophysiology experiment. *Statist. Med.* **19** 2067–2075.

Hu, F. and Rosenberger, W. F. (2003). Optimality, variability, power: Evaluating response-adaptive randomization procedures for treatment comparisons. *J. Amer. Statist. Assoc.* **98** 671–678. MR2011680

Johnson, N. L. and Kotz, S. (1977). *Urn Models and Their Applications.* Wiley, New York. MR488211

Mahmoud, H. M. and Smythe, R. T. (1991). On the distribution of leaves in rooted subtree of recursive trees. *Ann. Appl. Probab.* **1** 406–418. MR1111525

Matthews, P. C. and Rosenberger, W. F. (1997). Variance in randomized play-the-winner clinical trials. *Statist. Probab. Lett.* **35** 193–207.

Rosenberger, W. F. (1996). New directions in adaptive designs. *Statist. Sci.* **11** 137–149.

Rosenberger, W. F. (2002). Randomized urn models and sequential design (with discussion). *Sequential Anal.* **21** 1–21. MR1903097

Smythe, R. T. (1996). Central limit theorems for urn models. *Stochastic Process. Appl.* **65** 115–137. MR1422883

Wei, L. J. (1979). The generalized Pólya's urn design for sequential medical trials. *Ann. Statist.* **7** 291–296.

Wei, L. J. and Durham, S. (1978). The randomized play-the-winner rule in medical trials. *J. Amer. Statist. Assoc.* **73** 840–843. MR514157

Zelen, M. (1969). Play the winner rule and the controlled clinical trial. *J. Amer. Statist. Assoc.* **64** 131–146. MR240938

College of Mathematics and Statistics
Northeast Normal University
Changchun
China
and
Department of Statistics
  and Applied Probability
National University of Singapore
Singapore

Department of Statistics
Halsey Hall
University of Virginia
Charlottesville, Virginia 22904-4135
USA
e-mail: fh6e@virginia.edu